\documentclass[11pt]{amsart}

\usepackage{amsmath,amssymb,amsthm,amsrefs,a4wide}
\usepackage{latexsym}
\usepackage{indentfirst}
\usepackage{graphicx}
\usepackage{placeins}
\usepackage{booktabs}
\usepackage{algorithm}
\usepackage{algorithmic}
\usepackage{multirow}
\usepackage{color}
\usepackage{xcolor}

\theoremstyle{plain}

\theoremstyle{definition}

\theoremstyle{remark}

\DeclareMathOperator{\argmin}{argmin}

\newcommand{\E}{\mathrm{E}}

\newcommand{\bbm}{\begin{bmatrix}}
\newcommand{\ebm}{\end{bmatrix}}

\newcommand{\R}{\mathrm{R}}

\newcommand{\T}{\top}

\newcommand{\trho}{\rho^*}
\newcommand{\tpi}{\pi^*}
\newcommand{\tv}{v^*}
\newcommand{\ta}{a^*}

\newcommand{\bs}{\bar{s}}
\newcommand{\bv}{\bar{v}}

\begin{document}

\title[Optimization formulations of MDP] {A note on optimization formulations of \\Markov decision
  processes}

\author[]{Lexing Ying}
\address[Lexing Ying]{Department of Mathematics, Stanford University, Stanford, CA 94305}
\email{lexing@stanford.edu}

\author[]{Yuhua Zhu}
\address[Yuhua Zhu]{Department of Mathematics, Stanford University, Stanford, CA 94305}
\email{yuhuazhu@stanford.edu}

\thanks{The work is partially supported by the U.S. Department of Energy, Office of Science, Office
  of Advanced Scientific Computing Research, Scientific Discovery through Advanced Computing
  (SciDAC) program and also by the National Science Foundation under award DMS-1818449.  }


%\subjclass[2010]{65Z05, 82B28, 82B80.}

\begin{abstract}
  This note summarizes the optimization formulations used in the study of Markov decision
  processes. We consider both the discounted and undiscounted processes under the standard and
  the entropy-regularized settings. For each setting, we first summarize the primal, dual, and
  primal-dual problems of the linear programming formulation. We then detail the connections between
  these problems and other formulations for Markov decision processes such as the Bellman equation
  and the policy gradient method.
\end{abstract}

\maketitle
\def\ds{\displaystyle}

\def\S{\mathcal{S}}
\def\A{\mathcal{A}}
\def\M{\mathcal{M}}
\def\P{\mathcal{P}}
\def\R{\mathbb{R}}
\def\E{\mathbb{E}}
\def\g{\gamma}
\def\l{\left}
\def\r{\right}
\def\ll{\left\lVert}
\def\rl{\right\rVert}

%================================================
\section{Introduction}\label{sec:intro}

% fixed point perspective
Most of the algorithms of Markov decision processes (MDPs) are derived from the fixed-point
iteration of the Bellman equation \cite{bellman1966dynamic}. Examples include value iteration
\cite{puterman2014markov,bertsekas2018abstract,bertsekas1996neuro}, policy iteration
\cite{bellman1966dynamic,howard1960dynamic}, temporal difference (TD) learning
\cite{sutton1988learning}, Q-learning \cite{watkins1989learning}, etc. The analyses of these
algorithms in the tabular case and linear function approximation case often leverage the
contraction property of the Bellman operator. In the past decade or so, nonlinear approximations
such as neural networks have become more popular. However, for nonlinear function
approximations, this contraction property no longer holds, often resulting in instability. Many
variants and modifications have been proposed to stabilize the training, e.g., DQN
\cite{mnih2013playing}, A3C \cite{mnih2016asynchronous}. However, theoretical guarantees for these
algorithms are still missing.

% optimization perspective
A second perspective of studying MDPs is based on optimization. For nonlinear approximations,
optimization formulations are often more convenient both for algorithmic design and mathematical
analysis as they guarantee convergence to at least local minimums. Therefore in recent years, more
attention has been given to the optimization framework. One major direction is based on linear
programming (LP) \cite{puterman1990markov} and some recent developments include
\cite{wang2020randomized,abbasi2014linear,tang2019doubly,chen2018scalable}. Another direction is the
Bellman residual minimization (BRM) \cite{baird1995residual}, which includes algorithms based on the
primal-dual form of the BRM \cite{sutton2009fast,bhatnagar2009convergent,dai2018sbeed},
stochastic compositional gradient (SCGD) methods based on two-scale separation
\cite{wang2017accelerating,wang2017stochastic}, algorithms based on the smoothness of the underlying
transition dynamics \cite{zhu2020borrowing,ZhuQ2020,lee2019stochastic}. The convergence properties
of these algorithms have been studied in
\cite{liu2015finite,mahadevan2014proximal,sutton2008convergent,neu2017unified,ye2011simplex}.

% our contributions
The contribution of this note is two-fold. First, we summarize the LP problems used in the study of
MDPs. Many results in this note are well-known, but we were not able to find a place where
these results are summarized in a uniform framework.  Second, we point out the connections between
the LP problems and other MDP formulations, including the equivalence between the dual problem and the
policy gradient method and the equivalence between the primal problem and the Bellman equation.

%% some seemingly different formulations, including the
%% equivalence between the Bellman equation and the primal linear-programming formulation as
%% well as the equivalence between the policy gradient method and the dual linear-programming
%% formulation.

\subsection{Notation}

% introduce MDP's five components
A Markov decision process $\M$ with discrete state and action spaces is characterized by $\M =
(\S,\A,P,r,\g)$. Here $\S$ is the discrete state space, with each state usually denoted by $s$. $\A$
is the discrete action space, with each action usually denoted by $a\in\A$. Throughout the note,
$|\S|$ and $|\A|$ are used to denote the size of $\S$ and $\A$, respectively.  $P$ is a
third-order tensor where, for each action $a\in\A$, $P^a\in\R^{|\S|\times|\S|}$ is the transition
matrix between the states, i.e., $P^a_{st}$ is the probability of arriving at state $t$ if action
$a$ is taken at state $s$. $r$ is a second-order tensor where, for each action $a\in\A$, $r^a_s$ is
the reward at state $s$ if action $a$ is taken. Finally, $\g\in[0,1]$ is the discount factor.

Let $\Delta$ be the probability simplex over the space of actions, i.e., 
\[
\Delta = \left\{\eta=(\eta^a)_{a\in\A}: \sum_{a\in\A} \eta^a = 1\text{ and } \eta^a \geq 0\text{
  for }\forall a\in\A\right\}.
\]
The set of all valid policies is defined to be 
\[
\Delta^{|\S|} = \left\{\pi=(\pi_s)_{s\in\S}: \pi_s \in \Delta \text{ for }\forall s\in\S \right\}.
\]
For a policy $\pi \in \Delta^{|\S|}$, the transition matrix $P^\pi\in\R^{|\S|\times|\S|}$ under the
policy $\pi$ is defined as
\begin{equation}\label{def of Ppi}
  P^\pi_{st} = \sum_{a\in\A}P^a_{st}\pi^a_s,
\end{equation}
i.e., $P^\pi_{st}$ is the probability of arriving at state $t$ from state $s$ if policy
$\pi$ is taken. The reward $r^\pi\in\R^{|\S|}$ under the policy $\pi$ is
\begin{equation}\label{def of rpi}
  r^\pi_{s} = \sum_{a\in\A}r^a_s \pi^a_s,
\end{equation}
i.e., $r^\pi_s$ represents the expected reward at state $s$ if policy $\pi$ is taken. 

% discounted vs. undiscounted
%\LY{Now introduce the discounted and undiscounted MDPs.}

Each policy $\pi$ induces a discrete Markov process, where at each round $m$, an action $a_m$ is chosen at
state $s_m$ according to a particular policy $\pi$, and then the agent arrives at state $s_{m+1}$
according to the distribution of the transition matrix $P^{a_m}$ and receive a reward
$r^{a_m}_{s_m}$. The goal of an MDP problem is to maximize the cumulative reward among all possible
policies. Depending on whether $\g$ is strictly less than one or not, an MDP can either be
{\bf discounted} ($\g<1$) or {\bf undiscounted} ($\g=1$).
%% \begin{itemize}
%% \item The discounted MDP refers to the case $\g<1$. The reward occurs $t$ steps in the future from
%%   the current state is multiplied by $\g^t$ to indicate its importance to the current state. As $\g$
%%   approaches $1$, more future rewards will be taken into account.
%% \item The undiscounted setting refers to $\g = 1$. Taking into account the effect of each reward
%%   equally on the future evolution of the MDP.
%% \end{itemize}

%\LY{Now introduce the standard and regularized controls.}

When solving the MDPs, entropy regularizer has been proved quite useful in terms of exploration and
convergence \cite{peters2010relative,fox2015taming,schulman2015trust,mnih2016asynchronous}. In this
note, we adopt the following {\em negative conditional entropy} defined for non-negative
$\rho\in\R^{|\A|}$:
\begin{equation}\label{def of h}
  h(\rho) := \sum_{a\in\A} \rho^a \log \frac{\rho^a}{ \sum_{b\in\A} \rho^b}.
\end{equation}
This entropy $h(\rho)$ is both convex and homogeneous of degree one in $\rho$ (see for example
Appendix A.1 of \cite{neu2017unified}). This regularizer has been widely used in the literature
\cite{peters2010relative,fox2015taming,schulman2015trust,mnih2016asynchronous,dai2018sbeed,haarnoja2018soft},
Depending on whether this regularizer is used, we call an MDP either {\bf standard} or {\bf
  regularized}.

%% \begin{itemize}
%% \item The standard MDP refers to the usual one with no constraints or penalty on the reward.
%% \item The regularized MDP refers to the MDP with a regularized reward, where entropy is added to
%%   the original reward, or equivalently, is added to the objective function. 
%% \end{itemize}

\subsection{Outline}

The rest of the note is organized as follows. In Section \ref{sec:gn}, we first derive the primal,
dual, and primal-dual problems for the discounted standard MDP. We then show the equivalence between
the policy gradient algorithm and the dual problem as well as the equivalence between the Bellman
equation and the primal problem. Sections \ref{sec:gr}, \ref{sec:1n}, and \ref{sec:1r} address the
discounted regularized MDP, the undiscounted standard MDP, and the undiscounted regularized MDP,
respectively, by following the same outline.

%The connections to the corresponding Bellman equation and policy gradient in these MDP and
%the linear programming is also shown.

%================================================
\section{Discounted standard MDP}\label{sec:gn}

The discounted standard MDP is probably the most studied case in literature
\cite{puterman2014markov, sutton2018reinforcement}. For $\g \in (0,1)$, the value function under
policy $\pi$ is a vector $v^\pi \in\R^{|\S|}$, where $v^\pi_s$ represents the expected discounted cumulative
reward starting from state $s$ under the policy $\pi$, i.e.,
\[
v^\pi_s  = \E \l[\sum_{m=0}^\infty \g^m r^{a_m}_{s_m} | s_0 = s\r],
\]
where the expectation is taken over ${a_m\sim\pi_{s_m},\\ s_{m+1}\sim P^{a_m}_{s_m,\cdot}}$, for all
$m\geq 0$. The value function naturally satisfies the Bellman equation for any $s\in \S$:
\begin{equation*}
  v^\pi_s = r^\pi_s + \g \E^\pi [v^\pi_{s_1}|s_0 = s] = r^\pi_s + \g\sum_{t\in\S}P^\pi_{st}v^\pi_t.
\end{equation*}
The goal of an MDP problem is to find the maximum value function among all possible policies. 

%There are two equivalent perspectives to this task. One aims to find the maximum value function
%directly, which we refer as the primal problem. Another is to find the policy that maximizes the
%value function, which is called the dual problem.

%=======
\subsection{LP problems}\hfill
%\textcolor{blue}{(Shall we make every formulation in vector form?)}  $\textcolor{blue}{=> \tv =
%\max v}$ To solve the Bellman equation \eqref{eqn: opt_be} is equivalent to sovle
%\textcolor{blue}{(more details?)}

\begin{figure}[h]
  \includegraphics[width=0.7\linewidth]{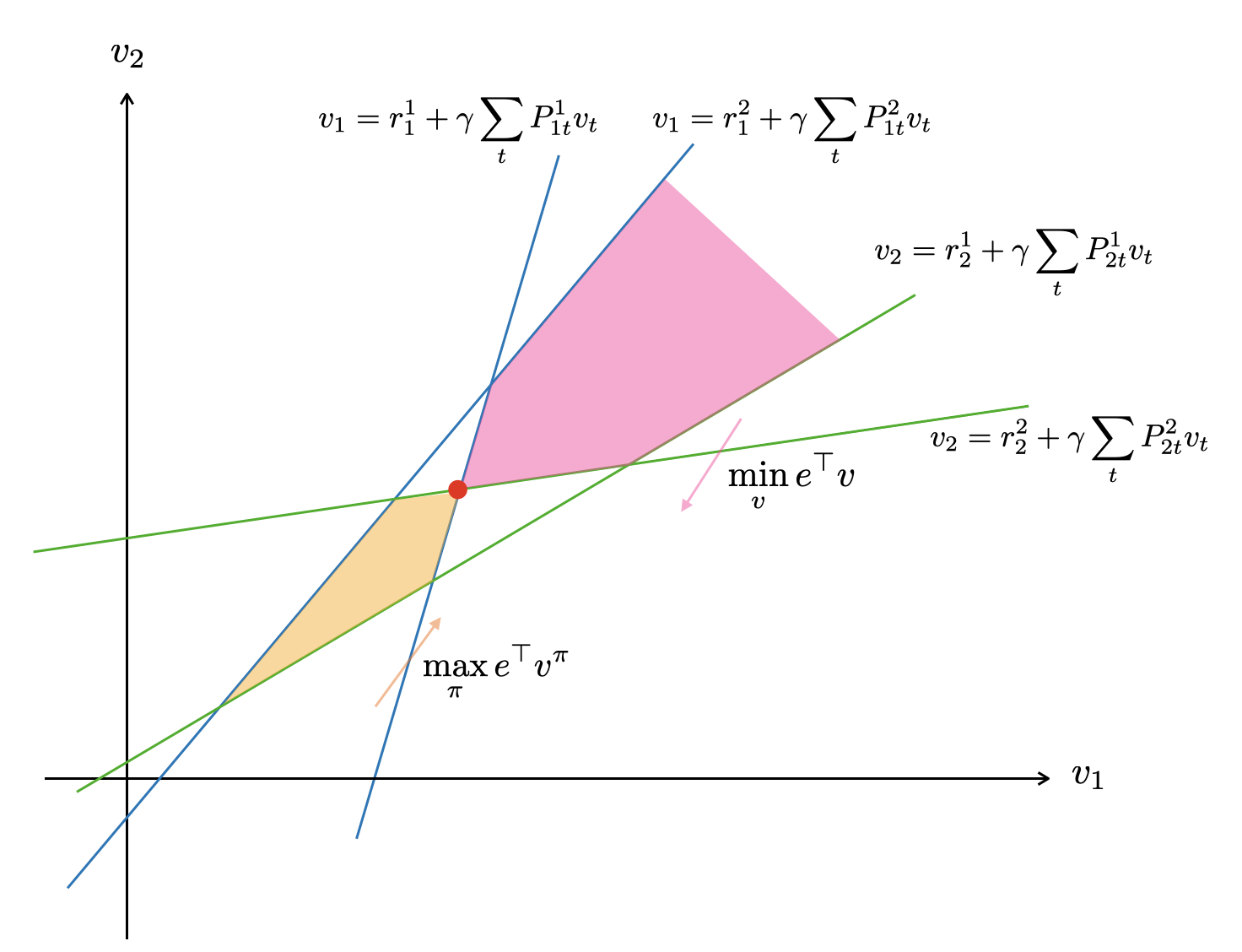}
  \caption{This plot interprets the primal and dual problem for an MDP with $2$ states and $2$
    actions. The pink region represents the constraints for the primal problem, while the yellow region
    represents the constraints for the dual problem. The pink and yellow arrows indicate the minimization
    and maximization directions, respectively. Both formulations end up at the same red point, but
    from the opposite directions. }
  \label{fig:itpt1}
\end{figure}

{\bf Primal problem.}  The primal problem of finding the maximum value function reads
\begin{equation}\label{eqn:gn_p}
  \min_{v} \  \sum_{s\in\S} e_s v_s,
  \text{ s.t. }
  \forall a, \forall s,\ r^a_s + \g \sum_{t\in\S} P^a_{st} v_t - v_s \le 0,
\end{equation}
where $e\in\R^{|\S|}$ is an arbitrary vector with positive entries.  An example with $2$ states and
$2$ actions is illustrated in Figure \ref{fig:itpt1}, where the pink region represents the
constraints and the pink arrow points to the minimization direction. The optimal solution of this
minimization problem is the red point. As shall see later in Section \ref{sec: equiv_gn}, the
optimal solution of the dual problem is the same red point, coming from the opposite direction.

{\bf Primal-dual problem.}  By introducing the Lagrangian multiplier $\mu^a_s$ for the inequality
constraints, we arrive at the primal-dual
\begin{equation}\label{eqn:gn_pd}
  \min_{v_s}\max_{\mu^a_s \ge 0}\  \sum_{s\in\S} e_s v_s + \sum_{s,a} (r^a_s+\g\sum_{t\in\S} P^a_{st}v_t-v_s)\mu^a_s,
\end{equation}
or equivalently in the matrix-vector notation
\[
\min_{v} \max_{\mu^a \ge 0}\  e^\T v + \sum_{a\in\A}(\mu^a)^\T (r^a + \g P^av - v),
\]
where $(\cdot)^\T$ stands for transpose. This primal-dual problem is, for example, used in
\cite{wang2020randomized}.

{\bf Dual problem.}  Since the minimum of $v_s$ is taken over a convex function and the maximum
of $\mu_s^a$ is over a concave function, one can exchange the order of minimum and maximum because
of the minimax theorem \cite{neumann1928theorie}. The primal-dual problem can also be written as
\[
\max_{\mu^a \ge 0} \min_{v} \  e^\T v + \sum_{a\in\A}(\mu^a)^\T (r^a + \g P^av - v).
\]
Taking derivative with respect to $v$ and setting it to be zero gives rise to
\[
e = -\sum_{a\in\A} (\g(P^a)^\T - I)\mu^a, \quad i.e., \quad \sum_{a\in\A} (I-\g(P^a)^\T)\mu^a = e.
\]
Hence the dual problem is
\begin{equation}\label{eqn:gn_d}
  \max_{\mu^a \ge 0}  \  \sum_{a\in\A} (r^a)^\T \mu^a,
  \text{ s.t. }
  \sum_{a\in\A} (I-\g(P^a)^\T)\mu^a = e.
\end{equation}
This dual problem is mentioned, for example, in  \cite{ye2011simplex}.

%=======
\subsection{Equivalences}\hfill
\label{sec: equiv_gn}

{\bf Dual problem and policy gradient.}  The dual problem \eqref{eqn:gn_d} is equivalent to the
policy gradient method. To see this, let us parameterize $\mu^a_s = w_s \pi^a_s$ with
$w_s=\sum_{a\in\A}\mu^a_s$. This ensures that $\pi\in\Delta^{|\S|}$ because
$\sum_{a\in\A}\pi^a_s=\sum_{a\in\A}\frac{\mu^a_s}{w_s} = 1$ and $\pi^a_s\geq 0$. By this new parameterization,
the constraints in dual become
\[
(I - \g (P^\pi)^\T) w = e,\quad\text{or}\quad
w = (I - \g (P^\pi)^\T)^{-1} e,
\]
where $P^\pi_{st}$ is defined in \eqref{def of Ppi} as the transition matrix under policy $\pi$.  By
denoting this $w$ as $w^\pi$ to indicate its $\pi$ dependence, we can write $\sum_{a\in\A} (r^a)^\T \mu^a =
(r^\pi)^\T w^\pi$. As a result, the dual problem \eqref{eqn:gn_d} can be written as
\[
  \max_{\pi\in\Delta^{|S|}} r^\pi \cdot (I - \g (P^\pi)^\T)^{-1} e,\quad\text{or}\quad
  \max_{\pi\in\Delta^{|S|}} e^\T (I - \g P^\pi)^{-1} r^\pi.
\]
It is clearly equivalent to the policy gradient method
\begin{equation}\label{eqn:gn_pg}
  \max_{\pi\in\Delta^{|\S|}} \ e^\T v^\pi, \text{ s.t. } v^\pi = r^\pi + \g P^\pi v^\pi,
\end{equation}
where we recall that $e\in\R^{|\S|}$ is any vector with positive entries. Therefore, {\bf the policy
  gradient method can be viewed as a nonlinear reparameterization of the dual LP problem}. This
understanding is also illustrated in Figure \ref{fig:itpt1}, where the yellow region represents the
constraints and the yellow arrow points to the maximum direction. Notice that both the primal and
dual problems end up at the same red point from opposite directions.

%The above formulation is the objective function for the policy gradient and we will show that the
%policy gradient objective is actually the dual problem of LP in a special parameterization.

%Therefore, {\bf the policy gradient objective is the dual problem in a special
%  parameterization.}

{\bf Primal problem and Bellman equation.} Next, we show that the primal problem
is equivalent to solving the Bellman equation %\YZ{(Shall we call it Bellman equation or Bellman equation?)}
\begin{equation}\label{eqn: opt_be}
  v_s =\max_{a\in\A} \left(r^a_s + \g\sum_{t\in\S}P^a_{st}v_t\right).
\end{equation}

{\em Bellman equation to primal problem.} The derivation from \eqref{eqn: opt_be} to \eqref{eqn:gn_p} can be found for example in
\cite{puterman2014markov}. We provide a short derivation here for completeness. Let $\tv$ be the solution to \eqref{eqn: opt_be}, then for
each $s$, there exists $\ta_s\in\A$ s.t.,
%\begin{equation}\label{eqn: gn_equiv1}
\[    \tv = r^{\ta} + \g P^{\ta}\tv. \]
%\end{equation} 
where $r^{\ta}_s \equiv r^{\ta_s}_s, P^{\ta}_{st} \equiv P^{\ta_s}_{st}$.  For any $v$ that
satisfies the constraints in the primal problem \eqref{eqn:gn_p}, the following inequality holds
%\begin{equation}\label{eqn: gn_equiv2}
\[  v \geq r^{\ta} + \g P^{\ta}v.\]
%\end{equation} 
Subtracting these two equations gives
\[
v-\tv \geq \g P^{\ta}(v - \tv).
\]
Since $P^{\ta}$ is a probability transition matrix, by maximum principle, one has $v - \tv \geq 0$, and thus 
\[
e^\T v \geq e^\T \tv
\]
for all $v$ satisfying the constraints in \eqref{eqn:gn_p}. This proves that $\tv$ is the minimizer
of the primal problem \eqref{eqn:gn_p}.

{\em Primal problem to Bellman equation.}  Let $\tv$ be the minimizer of the primal problem
\eqref{eqn:gn_p}. The KKT conditions for \eqref{eqn:gn_p} read
\[ %begin{equation}
  \l\{
  \begin{aligned}
    &r^a_s + \g\sum_{t\in\S}P^a_{st}\tv_t \leq \tv_s, \quad \text{for } \forall s, a;\\
    &\sum_{a\in\A}(\mu^a_s - \g\sum_{t\in\S}P^a_{ts}\mu^a_t) = e_s, \quad \text{for } \forall s;\\
    &\mu^a_s(r^a_s + \g\sum_{t\in\S}P^a_{st}\tv_t - \tv_s) = 0\quad \text{for } \forall s, a.
  \end{aligned}
  \r.
\] %end{equation}
First, we claim it is impossible that there exists $s$ such that for $\forall a$, $r^a_s +
\g\sum_{t\in\S}P^a_{st}\tv_t-\tv_s>0$. If it were true, then from the last equation, one would have
$\mu^a_s=0$ for all $a$. Let $\mu^a_s = w_s\pi^a_s$, then $w_s = 0$ for this $s$.  Inserting it into
the second equation, one has
\[
w_s - \g\sum_{t\in\S}P^\pi_{ts}w_t = - \g\sum_{t\in\S}P^\pi_{ts}w_t = e_s.
\]
Since both $P^\pi_{ts}, w_t \geq 0$ for all $t$, then LHS $\le 0$. However, as the RHS $e_s>0$, we
reach a contradiction. Therefore, the claim is true, i.e., There does not exist any $s$ such
that $\forall a$, $r^a_s+g\sum_{t\in\S}P^a_{st}\tv_t-\tv_s>0$.\\
Therefore for any fixed $s$, there exists $\ta_s$, s.t.,
\[
r^{\ta_s}_s + \g\sum_{t\in\S}P^{\ta_s}_{st}\tv_t = \tv_s.
\]
For all $a \neq \ta_s$, by the first equation, one has
\[
r^a_s + \g\sum_{t\in\S}P^a_{st}\tv_t \leq \tv_s.
\]
Combining the above two equations leads to
%\[
%\tv_s = r^{\ta_s}_s + \g\sum_{t\in\S}P^{\ta_s}_{st}\tv_t \geq r^a_s + \g\sum_{t\in\S}P^a_{st}\tv_t,
%\]
%which implies that
$\tv_s = \max_a(r^a_s + \g\sum_{t\in\S}P^a_{st}\tv_t)$. Therefore, the minimizer $\tv$ also satisfies the
Bellman equation \eqref{eqn: opt_be} for all $s$.

%\textcolor{gray}{It is easy to see, if $v \geq r + \g P^a v$, then $v\geq (I-\g P^a)^{-1}r$, so $v \geq v^a$ for $v^a := r + \g P^a v^a$ because $v^a = (I-\g P^a)^{-1}r$. Therefore the minimizer for \[\min v_s, \text{ s.t. } v \geq r + \g P^a v\] is $v^a_s$. This implies that the minimizer for \[\min v_s, \text{ s.t. } \forall a,\quad v \geq r + \g P^a v\]is $\tv_s = v^{a*}_s$, where $\ta = \text{argmax}_a v^a_s$. Because if $\tv_s = v^a_s$, for $v^a_s < v^{a*}_s$, then the constraint $v^a_s \geq r_s + \g\sum_{t\in\S}P^{a*}_{st}v_t$ gives $v^a_s \geq v^{a*}_s$, which contradicts with $v^a_s < v^{a*}_s$, so $\tv_s = \max_a v^a_s$. The objective function for \eqref{eqn:gn_p} is $e^\T v$, where $e_s>0$ for all $s$. Thererfore, $\tv_s = \max_a v^a_s$, or equivalently, $\tv_s = \max_a(r_s + \g \sum_{t\in\S}P^a_{st}\tv_t)$ should hold for all $s$, which means that the minimizer for \eqref{eqn:gn_p} satisfies the Bellman equation \eqref{eqn: opt_be}. }

%This is also why it is called primal Bellman equaiton as well. 

%================================================
\section{Discounted regularized MDP}
\label{sec:gr}

The discounted regularized MDP includes the negative conditional entropy
\[
h(\mu_s) = \sum_{a\in\A} \mu^a_s \log \frac{\mu^a_s}{\sum_{b\in\A} \mu^b_s}
\]
for each $s\in\S$ in the objective function.

%Note that there are also other form of regularizer such as relative entropy, which has been used in
%\cite{neu2017unified,belousov2017f,wu2019behavior,haarnoja2017reinforcement,geist2019theory}.

%=======
\subsection{LP problems}\hfill

{\bf Primal-dual problem.}  Let us introduce first in the primal-dual problem:
\begin{equation}\label{eqn:gr_pd}
  \min_{v_s} \sup_{\mu^a_s > 0} \sum_{s\in\S} e_s v_s + \sum_{s,a} \left(r^a_s+\g\sum_{t\in\S}
  P^a_{st}v_t-v_s\right)\mu^a_s - \sum_{s\in\S} h(\mu_s).
\end{equation}
As $h(\mu_s)$ is convex in $\mu_s$, this objective function is concave in $\mu^a_s$.  In connection
with the primal-dual problem \eqref{eqn:gn_pd} of the standard case, including the extra
entropic term allows for replacing the condition $\mu^a_s \ge 0$ with $\mu^a_s>0$. In the
literature, it is common for the entropy term to have a prefactor $\eta>0$. Here we simply assume
$\eta=1$ as one can always reduce to this case by rescaling the rewards $r^a_s$.

%=======
{\bf Primal problem.}  By introducing $\mu^a_s = w_s \pi^a_s$ with $w_s = \sum_{a\in\A} \mu^a_s$ and
$\pi\in\Delta^{|\S|}$, \eqref{eqn:gr_pd} becomes
\[ 
  \min_{v_s} \sup_{\pi \in \Delta^{|\S|},\pi^a_s,w_s > 0} \sum_{s\in\S} e_s v_s + \sum_{s,a} (r^a_s + \g\sum_{t\in\S}
  P^a_{st}v_t-v_s)\pi^a_sw_s - \sum_{s\in\S} w_s h(\pi_s),
\]
where we use the fact that $h(\cdot)$ is homogeneous of degree one. This is
equivalent to
\[
\min_{v_s} \left( \sum_{s\in\S} e_s v_s + \sup_{w_s> 0} \sum_{s\in\S} w_s \cdot \max_{\pi_s\in \Delta} \left(
\sum_{a\in\A} (r^a_s+\g\sum_{t\in\S} P^a_{st}v_t-v_s)\pi^a_s - h(\pi_s)\right) \right).
\]
Since the inner optimal over $\pi_s$ cannot lie on the boundary, it is the same to write the optimal
as $\max_{\pi_s\in\Delta}$ and $\sup_{\pi_s\in\Delta,\pi_s>0}$.  The primal problem of the above
minimax problem is then given by
\[
  \min_{v_s} \ e^\T v, \text{ s.t. } \forall s, \ \max_{\pi_s\in \Delta} \sum_{a\in\A} (r^a_s+\g\sum_{t\in\S}
  P^a_{st}v_t-v_s)\pi^a_s - h(\pi_s) \leq 0,
\]
%The above formulation is also the LP formulation for the regularized Bellman equation
%\eqref{eqn:gr_be} (\textcolor{blue}{More details?})
or equivalently, 
\begin{equation}\label{eqn:gr_p0}
    \min_{v}e^\T v, \text{ s.t. } \max_{\pi} r^\pi + \g P^\pi v - h^\pi \leq v.
\end{equation}
Note that the maximization in the constraint
\[
\max_{\pi_s\in \Delta} \left( \sum_{a\in\A} (r^a_s+\g\sum_{t\in\S} P^a_{st}v_t-v_s)\pi^a_s - h(\pi_s)\right)
\]
is in the form of the Gibbs variational principle. Therefore, the optimizer for this maximization is
\[
\pi^a_s = \frac{\exp(r^a_s+\g\sum_{t\in\S} P^a_{st}v_t-v_s)}{Z_s}
\]
where $Z_s = \sum_{a\in\A} \exp(r^a_s+\g\sum_{t\in\S} P^a_{st}v_t-v_s)$ is the normalization factor and the
maximal value is $\log Z_s$. Hence, the constraint is equivalent to $\log Z_s \leq 0$, i.e.,
%d $w_s\ge 0$, the maximization will give infinity unless the optimal value $\log Z_s \le 0$, i.e.,
%$Z_s \le 1$. This is
\[
\sum_{a\in\A}  \exp(r^a_s+\g\sum_{t\in\S} P^a_{st}v_t-v_s) \le 1
\]
i.e.
\[
e^{-v_s} \cdot \sum_{a\in\A}  \exp(r^a_s+\g\sum_{t\in\S} P^a_{st}v_t) \le 1,\quad\text{or}\quad
\sum_{a\in\A}  \exp(r^a_s+\g\sum_{t\in\S} P^a_{st}v_t) \le e^{v_s}.
\]
Taking log on both sides leads to 
\[
v_s \ge \log\left(\sum_{a\in\A}  \exp(r^a_s+\g\sum_{t\in\S} P^a_{st}v_t)\right).
\]
This implies that the primal problem \eqref{eqn:gr_p0} is equivalent to
\begin{equation}\label{eqn:gr_p}
  \min_{v_s}\ e^\T v, \text{ s.t. } \forall s, v_s \ge \log\left(\sum_{a\in\A} \exp(r^a_s+\g\sum_{t\in\S}
  P^a_{st}v_t)\right).
\end{equation}

%{\color{blue} If (10) becomes "=0", then here it should be "=" as well. Besides, the operator is an
%constractive operator, so it only has one fixed point, which means the below equation optimization
%can be direction written as the constraint = 0.}

%The above formulation is also the LP formulation for the log-sum-exp form of the optimal value
%function \eqref{eqn:gr_logsum}. (\textcolor{blue}{More details?})Therefore, \eqref{eqn:gr_be} is
%equivalent to \eqref{eqn:gr_logsum}, which is also used in \cite{dai2018sbeed}.

%=======
{\bf Dual problem.}  The supremum and minimum in the primal-dual problem \eqref{eqn:gr_pd}
can be exchanged because the objective function is convex in $v$ and concave in $\mu$.  Now one has,
\[
  \sup_{\mu^a_s > 0} \min_{v_s}  \sum_{s\in\S} e_s v_s + \sum_{s,a} (r_s+\g\sum_{t\in\S} P^a_{st}v_t-v_s)\mu^a_s - \sum_{s\in\S} h(\mu_s).
\]
Taking derivative in $v$ gives $\sum_{a\in\A}(I-\g (P^a)^\T)\mu^a = e$. Therefore, the dual problem
takes the form
\begin{equation}\label{eqn:gr_d}
  \sup_{\mu^a>0} \sum_{a\in\A} (r^a)^\T \mu^a - \sum_{s\in\S} h(\mu_s), \text{ s.t. } \sum_{a\in\A}(I-\g (P^a)^\T)\mu^a
  = e.
\end{equation}

%=======
\subsection{Equivalences}\hfill
\label{sec: equiv_gr}

{\bf Dual problem and policy gradient.}  We claim that the dual problem \eqref{eqn:gr_d} is again
equivalent to the policy gradient method.  As before, let us parameterize $\mu^a_s = w_s \pi^a_s$
with $w_s = \sum_{a\in\A}\mu^a_s$ and $\pi\in\Delta^{|\S|}$. Then the constraints in \eqref{eqn:gr_d}
become
\[
\forall s, \ \sum_{a\in\A}\pi^a_sw_s - \g\sum_{a,t}P^a_{ts}\pi^a_tw_t = e_s, \quad\text{or}\quad w = (I-\g(P^\pi)^\T)^{-1}e.
\]
By denoting the solution $w$ by $w^\pi$ to show its $\pi$ dependence, one can rewrite
\eqref{eqn:gr_d} as
%\begin{equation}\label{eqn:gr_d}\sup_{\mu^a_s > 0}  \sum_{s\in\S} r_s (\sum_{a\in\A}\mu^a_s)- \sum_{s\in\S} h(\mu_s)\text{ s.t. }(I-\g(P^\pi)^\T)w = e.\end{equation}
%\[
%\max_{\pi^a_s,w_s} \min_{v_s} \sum_{s\in\S} e_s v_s + \sum_{s,a} (r_s+\g\sum_{t\in\S} P^a_{st}v_t-v_s)\pi^a_sw_s - \sum_{s\in\S} w_sh(\pi_s).
%\]
%Taking derivative in $v_s$ leads to
%In terms of $\mu^a_s =\pi^a_s w_s$, we have\[w[\pi] = (I - \g (P^\pi)^\T)^{-1} e\] as before. Plugging it in gives 
\[
\max_{\pi\in \Delta^{|\S|}} r^\pi \cdot w^\pi - \sum_{s\in\S} w^\pi_s \left(\sum_{a\in\A} \pi^a_s \log \pi^a_s\right)
\]
By further introducing $h^\pi \in\R^{|\S|}$ as the vector with entry $h^\pi_s = h(\pi_s) = \sum_{a\in\A}
\pi^a_s \log \pi^a_s$, we transform \eqref{eqn:gr_d} to
\[ %begin{equation}\label{eqn:gr_pg}
  \max_{\pi\in \Delta^{|\S|}} \left( r^\pi-h^\pi \right) \cdot w^\pi, \quad\text{or}\quad
  \max_{\pi\in \Delta^{|\S|}} e^\T (I - \g P^\pi)^{-1} \left( r^\pi-h^\pi \right).
\] %end{equation}

%% %and introduce 
%% Writing out $w[\pi]$ explicitly gives a equivalent dual problem,
%% \begin{equation}\label{eqn:gr_pg}
%%   \max_{\pi\in \Delta^{|\S|}} e^\T (I - \g P^\pi)^{-1} \left( r-h(\pi) \right),
%% \end{equation}
%% where $h(\mu)\in\R^s$ is a vector with the $s$-th element $h(\mu_s)$.

We can also view $r^\pi - h^\pi$ as a regularized reward by subtracting the entropy function
$h^\pi$. The value function $v^\pi = \E[\sum_{m\geq 0} \g^m(r^{a_m}_{s_m} - h(\pi_{s_m}))]$ under
the new reward satisfies the regularized Bellman equation $ v^\pi = r^\pi-h^\pi + \g P^\pi
v^\pi$. Hence, the policy gradient method of this regularized discounted MDP is
\begin{equation}\label{eqn:gr_pgpg}
  \max_{\pi\in \Delta^{|\S|}} e^\T v^\pi, \quad s.t., \quad v^\pi = r^\pi - h^\pi + \g P^\pi v^\pi,
\end{equation}
which is clearly equivalent.
%\eqref{eqn:gr_p}  and \eqref{eqn:gr_pg}. 

{\bf Primal problem and Bellman equation.}  The regularized Bellman equation is 
\begin{equation}\label{eqn:gr_be}
  v = \max_{\pi\in \Delta^{|\S|}} r^\pi + \g  P^\pi v - h^\pi.
\end{equation} 
In each component, 
\[
v_s = \max_{\pi_s\in \Delta} \sum_{a\in\A}(r^a_s + \g \sum_{t\in\S}P^a_{st} v_t))\pi^a_s - h(\pi_s).
\]
By the Gibbs variational principle, the RHS is equal to
$\log(\exp(\sum_{a\in\A}(r_s+\g\sum_{t\in\S}P^a_{st}v_t)))$. Therefore, the regularized Bellman equation could
also be written as the following log-sum-exp form
\begin{equation}\label{eqn:gr_logsum}   
  v_s = \log\left(\sum_{a\in\A}  \exp\left(r^a_s+\g\sum_{t\in\S} P^a_{st}v_t\right)\right).
\end{equation} 
Below we show that the primal problem \eqref{eqn:gr_p0} is equivalent to solving \eqref{eqn:gr_be}.

%% Besides, the maximum in \eqref{eqn:gr_be} could be written as $v_s^* = \max_{\pi_s\in \Delta} \sum_{a\in\A}(r^a_s + \g \sum_{t\in\S}P^a_{st} v_t))\pi^a_s - h(\pi_s)$, so by the Gibbs variational
%% principle, the maximum  $v_s^*$ could also be written in form of $v_s^* = \log(\exp(\sum_{a\in\A}(r_s + \g \sum_{t\in\S}P^a_{st} v_t)))$. Therefore, the optimal value function could also be written as the log-sum-exp form, 
%% \begin{equation}\label{eqn:gr_logsum}   
%% v_s^* = \log\left(\sum_{a\in\A}  \exp(r_s+\g\sum_{t\in\S} P^a_{st}v_t^*)\right).
%% \end{equation} 
%% Hence the primal problem \eqref{eqn:gr_p0} is also equivalent to the above Bellman equation in log-sum-exp form.

{\em Bellman equation to primal problem.} Let $\tv$ be the solution of \eqref{eqn:gr_be}. Then there
exists $\tpi$ s.t.,
\[
  \tv = r^{\tpi} + \g P^{\tpi}\tv - h^{\tpi}.
\]
For any $v$ that satisfies the constraints of the primal problem, the following inequality holds
\[
  v \geq r^{\tpi} + \g P^{\tpi}v - h^{\tpi}.
\]
Subtracting these two equations gives rise to
\[
v-\tv \geq \g P^{\tpi}(v - \tv).
\]
Again by maximum principle, one has $v - \tv \geq 0$ and hence
\[
e^\T v \geq e^\T \tv
\]
for all $v$ that satisfying the constraints in \eqref{eqn:gr_p0}. This proves that $\tv$ is the
minimizer of the primal problem \eqref{eqn:gr_p0}.

{\em Primal problem to Bellman equation.}  Let $\tv$ be the minimizer of the primal problem
\eqref{eqn:gr_p0}. We now prove $\tv$ is also the solution to the Bellman equation \eqref{eqn:gr_be}
by contradiction. Assume that $\tv$ does not satisfy \eqref{eqn:gr_be}. Then there must exist
$\bs$, s.t. for $\forall \pi$
\[ %begin{equation}\label{eqn: gr_equiv3}
    \tv_{\bs} \geq (r^\pi + \g P^\pi \tv - h^\pi)_{\bs} + \delta
\] %end{equation}
with some constant $\delta>0$. Let us define $\bv$ s.t., $\bv_{\bs} = \tv_{\bs} - \delta$ and $\bv_s
= \tv_s$ for $s\neq\bs$. We claim that for $\forall\pi$
\[
\bv_s \geq (r^\pi + \g P^\pi \bv - h^\pi)_s, \quad \forall s.
\]
First, for $s\neq \bs$ the above inequality holds because $\bv_s = \tv_s$ and $\tv_s$ satisfies the
constraints in the primal problem \eqref{eqn:gr_p0}.  For $s = \bs$, one has
$\bv_{\bs}=\tv_{\bs}-\delta \geq (r^\pi + \g P^\pi \tv- h^\pi)_{\bs}$. Since
\[
(P^\pi \tv)_{\bs} = \sum_{t\in\S}P^\pi_{\bs t}\tv_t = \sum_{t\neq\bs}P^\pi_{\bs t}\tv_t + P^\pi_{\bs
  \bs}(\tv_{\bs}- \delta) + \delta = \sum_{t\in\S} P^\pi_{\bs t}\bv_t + \delta,
\]
$\bv_{\bs} \geq (r^\pi + \g P^\pi \bv- h^\pi)_{\bs} + \delta$. This completes the proof of the
claim. This means that $\bv$ also satisfies the constraints of the primal problem, but $e^\T\bv <
e^\T\tv$ by construction, which contradicts with $\tv$ being the minimizer of the primal
problem. Therefore, the assumption is wrong and $\tv$ satisfies \eqref{eqn:gr_be}.

%The above equation is actually equivalent to the primal problem \eqref{eqn:gr_p0}. For all $v\geq\max_{\pi\in \Delta^{|\S|}} r^\pi + \g P^\pi v - h^\pi$ and $\tv=\max_{\pi\in \Delta^{|\S|}} r^\pi + \g P^\pi \tv - h^\pi$, one has $v \geq \tv$, therefore the minimizer of the primal problem $e^\T v$ is $\tv$. So the solution to \eqref{eqn:gr_p0} satisfies the optimal regularized Bellman equation \eqref{eqn:gr_be}. On the other hand, the solution to \eqref{eqn:gr_be} is also the minimizer of the regularized primal linear programming \eqref{eqn:gr_p0}.

%================================================
\section{Undiscounted standard MDP}
\label{sec:1n}
In this section, we consider the MDP without discounts, i.e., $\g = 1$. Besides, we assume the MDP
is unichain, i.e., for each policy $\pi$, the MDP induced by policy $\pi$ is ergodic
\cite{ortner2006short}. Let $\rho^\pi\in\R$ be the average reward under policy $\pi$,
\[
\rho^\pi = \lim_{T\to\infty} \E\left[ \frac{1}{T}\sum_{m=1}^T r^{a_m}_{s_m} \right],
\]
where the expectation is taken over $a_m\sim \pi_{s_m}, s_{m+1} \sim P^{a_m}_{s_m,\cdot}$. Notice
for an ergodic Markov process, the average-reward is the same for any
initial states because the stationary distribution is unique invariant, strictly positive and independent of initial states (Chapter 6 of \cite{kemeny2012denumerable}).  Let $w^\pi\in\R^{|\S|}>0$ be the stationary distribution induced by policy $\pi$
satisfying $w^\pi_s = \sum_{t\in\S} P^\pi_{ts}w^\pi_t$. 
%Note that for unichain MDP, each policy $\pi$ induces a unique stationary distribution no matter what the initial distribution is, and the stationary distribution is strictly positive. 
Then $w^\pi$ is the only solution to $(I - (P^\pi)^\T)w^\pi = 0$
with $\sum_{s\in\S}w^\pi_s = 1$.  The average reward under policy $\pi$ could also be written as
\[
\rho^\pi = (r^\pi)^\T w^\pi.
\]

%It also satisfies the following average-reward Bellman equaiton,
%\[
%V^\pi_s = r_s  - \rho^\pi + \g \sum_{t\in\S}P^\pi_{st}V_t.
%\]

Define the value function for the non-discount MDP as $v^\pi_s = \E_\pi[\sum_{m=1}^\infty
  (r^{a_m}_{s_m} - \rho^\pi) | s_0 = s]$, then the value function satisfies the average-reward
Bellman equation:
\begin{equation}\label{eqn:1n_be}
  v^\pi_s = r^\pi_s - \rho^\pi + \sum_{t\in\S}P^\pi_{st}v_t^\pi
\end{equation}
with $\sum_{s\in\S} v^\pi_sw^\pi_s = 1$. Note that without this constrain, there are infinitely many $v^\pi$,
e.g., $v^\pi + C$ for any constant $C$ still satisfies the above equation. The goal here is to find the
maximum average reward among all possible policies.

%=======
\subsection{LP problems}\hfill

{\bf Primal problem.}  The primal problem is
\begin{equation}\label{eqn:1n_p}
  \min_{v_s,\rho} \rho\ 
  \text{ s.t. }
  \forall a,\forall s, \  r^a_s + \sum_{t\in\S} P^a_{st} v_t - v_s - \rho \le 0.
\end{equation}
This can be found for example in  \cite{neu2017unified}.

%=======
{\bf Primal-dual problem.}
By including the Lagrangian multiplier $\mu^a_s$ for the inequality constraints, one obtains the
primal-dual problem
\begin{equation}\label{eqn:1n_pd}
  \min_{v_s,\rho} \max_{\mu^a_s \ge 0} \rho + \sum_{s,a} (r^a_s+\sum_{t\in\S} P^a_{st}v_t-v_s-\rho)\mu^a_s,
\end{equation}
or equivalently, in the matrix-vector notation
\[
\min_{v_s,\rho} \max_{\mu^a_s \ge 0}\  \rho + \sum_{a\in\A} (\mu^a)^\T(r^a+ P^av-v-\rho{\bf 1}),
\]
where ${\bf 1}$ is the $|S|$-dimensional vector with all elements equal to $1$.

%=======
{\bf Dual problem.} To get the dual problem, we take the derivative with respect to $\rho$ to get
\[
1 - \sum_{s,a} \mu^a_s = 0.
\]
Taking the derivative with respect to $v$ leads to
\[
\sum_{a\in\A} (I-(P^a)^\T) \mu^a = 0.
\]
Hence the dual problem is
\begin{equation}\label{eqn:1n_d}
  \max_{\mu^a_s \ge 0} \sum_{a\in\A} (r^a)^\T \mu^a, \text{ s.t. } \sum_{a\in\A} (I-(P^a)^\T)\mu^a = 0, \quad 1 -
  \sum_{s,a} \mu^a_s = 0.
\end{equation}

%=======
\subsection{Equivalences}\hfill

{\bf Dual problem and policy gradient.}  The dual problem \eqref{eqn:1n_d} is equivalent to the
policy gradient. Let us again parameterize $\mu^a_s = w_s \pi^a_s$ with $w_s = \sum_{a\in\A}\mu^a_s$ and
$\pi\in\Delta^{|\S|}$.  By the new parameterization, the constraints become
\[\sum_{a\in\A} (I-(P^a)^\T)\mu^a = 0 \Rightarrow (I - (P^\pi)^\T) w = 0.\]$1 - \sum_{s,a} \mu^a_s = 0$ also implies that $1-\sum_{s\in\S} w_s=0$. 
%\[(I - (P^\pi)^\T) w = 0, \quad 1-\sum_{s\in\S} w_s=0.\] 
Together we conclude that $w$ is the
stationary distribution induced by $\pi$. By denoting this $w$ by $w^\pi$, we can write the dual
problem as
\begin{equation}\label{eqn:1n_pg}
  \max_{\pi\in\Delta^{|\S|}}  r^\pi \cdot w^\pi,
\end{equation}
which is exactly the optimization formulation of the policy gradient method.

%% where the maximizer $\tpi$ is the best policy that maximizes $\rho^\pi = r^\pi \cdot w^\pi$, where
%% $w^\pi$ is the unique stationary distribution under transition matrix $P^\pi$.  
%% %finding the best policy is equivalent to finding the best stationary distribution,
%% %%\textcolor{blue}{(It seems one stationary distribution not necessarily comes from one unique
%% %policy? Does that mean the optimal policy is not unique?)}
%% On the other hand, the maximum average reward under the best policy could be solved by the following
%% policy gradient method,
%% \begin{equation}\label{eqn:ln_pgpg}
%%   \rho^* = \max_{\pi \in \Delta^{|\S|}} r^\pi\cdot w^\pi,
%% \end{equation}
%% which is equivalent to the dual problem of the averaged-reward LP \eqref{eqn:1n_pg}.

{\bf Primal problem and Bellman equation.} Next, we show that the primal problem \eqref{eqn:1n_p} is
equivalent to the average reward Bellman equation for $v,\rho$
\begin{equation}\label{eqn:ar_optbe}
  v_s = \max_a \left(r^a_s - \rho + \sum_{t\in\S}P^a_{st}v_t \right), \quad s\in\S.
\end{equation}
%with $\sum_{s\in\S} \tv_sw^*_s = 1$.

{\em Bellman equation to primal problem.} The derivation from \eqref{eqn:ar_optbe} to
\eqref{eqn:1n_p} can be found for example in \cite{puterman2014markov}. Here we provide a short
proof for completeness. Let $\tv,\trho$ be the solution to the Bellman equation
\eqref{eqn:ar_optbe}, then for all $s$, there exists $\ta_s$ s.t.,
\[ %begin{equation}\label{eqn: 1n_equiv1}
\tv = r^{\ta} - \trho{\bf 1} + P^{\ta}\tv,
\] %end{equation} 
where $r^{\ta}_s \equiv r_s^{\ta_s}, P^{\ta}_{st} \equiv P^{\ta_s}_{st}$.  For any $v,\rho$ that
satisfy the constraints in the primal problem \eqref{eqn:1n_p}, the following inequality holds
\[ %begin{equation}\label{eqn: 1n_equiv2}
v \geq r^{\ta} - \rho{\bf 1} + P^{\ta}v.
\] %end{equation} 
Subtracting these two equations gives
\[ %begin{equation}\label{eqn: 1n_equiv3}
  v-\tv \geq P^{\ta}(v - \tv) - (\rho - \trho){\bf 1}
\] %end{equation}
Let $w^*$ be the stationary distribution induced by the policy $\pi_s^a = \l\{\begin{aligned}
 1, a = \ta_s\\0, a\neq \ta_s
\end{aligned}\r.$, then $(w^*)^\T = (w^*)^\T P^{\ta}$. Multiplying $(w^*)^\T$ to the last equation yields, 
\[
(w^*)^\T(\rho - \trho){\bf 1} \geq 0
\]
Since we assume the MDP is unichain, the stationary distribution $w^*$ for any policy is strictly
positive. This implies that
\[
\rho  \geq \trho
\]
holds for all $v, \rho$ satisfying the constraints in \eqref{eqn:1n_p}. This proves that $\tv,\trho$
is the minimizer of the primal problem \eqref{eqn:1n_p}.

{\em Primal problem to Bellman equation.}  Let $(\tv,\trho)$ be the minimizer of the primal problem
\eqref{eqn:1n_p}. We now show that $(\tv,\trho)$ also satisfies the average reward bellman equation
\eqref{eqn:ar_optbe}. The KKT conditions of \eqref{eqn:1n_p} are
\begin{equation*}
  \l\{
  \begin{aligned}
    &r^a_s + \sum_{t\in\S}P^a_{st}\tv_t - \tv_s\leq \trho, \quad \text{for } \forall s, a;\\
    &1 - \sum_{s,a}\mu^a_s = 0;\\
    &\sum_{a\in\A}(\mu^a_s - \sum_{t\in\S}P^a_{ts}\mu^a_t) = 0, \quad \text{for } \forall s;\\
    &\mu^a_s(r^a_s + \sum_{t\in\S}P^a_{st}\tv_t - \tv_s - \trho) = 0\quad \text{for } \forall s, a.
  \end{aligned}
  \r.
\end{equation*}
First we claim that it is impossible that there exists $s$ s.t. $\forall a$, $r^a_s + \sum_{t\in\S}P^a_{st}\tv_t
- \tv_s-\trho < 0$. Let $\mu^a_s = w_s\pi^a_s$ with $w_s = \sum_{a\in\A} \mu^a_s$ and $\pi_s\in\Delta$ for
all $s$. Plugging it into the second and third equation gives,
\[
1 - \sum_{s\in\S}w_s = 0,\quad w - (P^\pi)^\T w = 0.
\]
Therefore $w$ is the unique stationary distribution induced by the policy $\pi$. Since we assume the
MDP is unichain, the stationary distribution is strictly positive. If there exists $s$, s.t. for
$\forall a$, $r^a_s + \sum_{t\in\S}P^a_{st}\tv_t - \tv_s - \trho < 0$, then by the last equation of the
KKT condition, $\mu^a_s = w_s\pi^a_s= 0$ implies $w_s = 0$ for this $s$, which contradicts with the
unichain assumption. Therefore, the claim is true, i.e., there does not exist any $s$ s.t. $\forall a$, $r^a_s+g\sum_{t\in\S}P^a_{st}\tv_t-\tv_s>0$.\\
Therefore, for $\forall s$, there always exists $\ta_s$, s.t.,
\[
r^{\ta_s}_s + \sum_{t\in\S}P^{\ta_s}_{st}\tv_t - \tv_s = \trho.
\] 
By the first equation, for all $a\neq \ta$,
\[
r^a_s + \sum_{t\in\S}P^{a}_{st}\tv_t - \tv_s  \leq  \trho.
\]
Combining the above two equations gives $r^a_s + \sum_{t\in\S}P^{a}_{st}\tv_t - \tv_s \leq r_s +
\sum_{t\in\S}P^{\ta_s}_{st}\tv_t - \tv_s= \trho.$ This is equivalent to $\trho = \max_a r^a_s +
\sum_{t\in\S}P^{a}_{st}\tv _t - \tv_s$, the average reward Bellman equation
\eqref{eqn:ar_optbe}.

%================================================
\section{Undiscounted regularized MDP}
\label{sec:1r}
%For a probability density $\mu_s$ over actions, we define the normalized version $\hat{\mu_s}$
%\[
%\hat{\mu}^a_s = \frac{\mu^a_s}{ \sum_b \mu^b_s}
%\]
%and introduce the homogeneous entropy function
%\[
%h(\mu_s) = \sum_{a\in\A} \mu^a_s \log \hat{\mu}^a_s.
%\]
%Notice that the first $\mu$ is the usual dual variable while the second $\hat{\mu}$ is the normalized version.

%=======
\subsection{LP problems}\hfill

{\bf Primal-dual problem.}  We again use the the negative conditional entropy
\[
h(\mu_s) = \sum_{a\in\A} \mu^a_s \log \frac{\mu^a_s}{\sum_b \mu^b_s}
\]
as the regularizer.  The primal-dual problem of the undiscounted regularized MDP is
\begin{equation}\label{eqn:1r_pd}
  \min_{v_s,\rho} \sup_{\mu^a_s>0} \rho + \sum_{s,a} (r^a_s+\sum_{t\in\S} P^a_{st}v_t-v_s-\rho)\mu^a_s -
  \sum_{s\in\S} h(\mu_s).
\end{equation}

%=======
{\bf Primal problem.}  By introducing $\mu^a_s = w_s \pi^a_s$, one obtains
\[
  \min_{v_s,\rho} \sup_{\pi\in\Delta^{|\S|},pi_s^a,w_s> 0} \rho + \sum_{s,a} (r^a_s+\sum_{t\in\S}
  P^a_{st}v_t-v_s-\rho)\pi^a_sw_s - \sum_{s\in\S} w_s h(\pi_s),
\]
which is equivalent to 
\[
\min_{v_s,\rho} \left( \rho + \sup_{w_s>0}
\sum_{s\in\S} w_s \cdot \max_{\pi_s\in\Delta}
\left( \sum_{a\in\A} (r^a_s+\sum_{t\in\S} P^a_{st}v_t-v_s-\rho)\pi^a_s - h(\pi_s)\right) \right).
\]
The primal problem of the above minimax problem is 
\[
  \min_{v_s,\rho} \rho \ \text{ s.t. } \forall s, \ \max_{\pi_s\in\Delta} \left( \sum_{a\in\A}
  (r^a_s+\sum_{t\in\S} P^a_{st}v_t-v_s-\rho)\pi^a_s - h(\pi_s)\right) \leq 0,
\]
or equivalently, 
\begin{equation}\label{eqn:1r_p0}
  \min_{v,\rho}\rho\  \text{ s.t. } \max_{\pi} r^\pi  - \rho{\bf 1} + P^\pi v - h^\pi \leq v.
\end{equation}
Note that the constraint 
\[\max_{\pi_s\in\Delta} \left( \sum_{a\in\A}
(r^a_s+\g\sum_{t\in\S} P^a_{st}v_t-v_s-\rho)\pi^a_s-h(\pi_s)\right)\] is in the form of the Gibbs
variational principle. The optimizer for this maximization is
\[
\pi^a_s = \frac{\exp(r^a_s+\sum_{t\in\S} P^a_{st}v_t-v_s-\rho)}{Z_s},
\]
where $Z_s$ is the normalization factor, the optimal value is $\log Z_s$.  
%However, due to the $w_s$ term in front of it and $w_s\ge 0$, the maximization will give infinity
%unless the optimal value
Therefore, the constraint is equivalent to $\log Z_s \le 0$, i.e.,
\[
\sum_{a\in\A} \exp(r^a_s+\sum_{t\in\S} P^a_{st}v_t-v_s-\rho) \le 1
\]
i.e.
\[
e^{-v_s} \cdot \sum_{a\in\A} \exp(r^a_s+\sum_{t\in\S} P^a_{st}v_t-\rho) \le 1,\quad\text{or}\quad \sum_{a\in\A}
\exp(r^a_s+\sum_{t\in\S} P^a_{st}v_t-\rho) \le e^{v_s}.
\]
Taking log gives
\[
v_s \ge \log\left(\sum_{a\in\A}  \exp(r^a_s+\sum_{t\in\S} P^a_{st}v_t-\rho)\right).
\]
Hence the primal problem can also be written as
\begin{equation}\label{eqn:1r_p}
  \min_{v_s,\rho} \rho\ 
  \text{ s.t. }
  v_s \ge \log\left(\sum_{a\in\A}  \exp(r^a_s+\sum_{t\in\S} P^a_{st}v_t-\rho)\right),
\end{equation}
which is an alternative formulation of the primal problem. Notice that a similar problem with
equality constraints was instead derived in \cite{neu2017unified}. In practice, the inequality
constraints as in \eqref{eqn:1r_p} are often preferred since the feasibility set is then convex.

%=======
{\bf Dual problem.}  The supremum and minimum in the primal-dual problem \eqref{eqn:gr_pd} can be
exchanged because the objective function is convex in $v,\rho$ and concave in $\mu$. Then one has,
\[
  \sup_{\mu^a_s>0} \min_{v_s,\rho}  \rho + \sum_{s,a} (r^a_s+\sum_{t\in\S} P^a_{st}v_t-v_s-\rho)\mu^a_s - \sum_{s\in\S} h(\mu_s).
\]
Taking derivatives in $\rho$ and $v$ leads to
\[
1-\sum_{s,a} \mu_s = 0, \quad \sum_{a\in\A}(I - (P^a)^\T)\mu^a = 0.
\]
Hence the dual problem in terms of $\mu^a_s$ is
\begin{equation}\label{eqn:1r_d}
  \sup_{\mu^a>0} \sum_{a\in\A} (r^a)^\T \mu^a- \sum_{s\in\S} h(\mu_s)\ \text{ s.t. }\sum_{a\in\A} (I-(P^a)^\T)\mu^a =
  0,\quad 1 - \sum_{s,a} \mu^a_s = 0.
\end{equation}

%Rewriting it in terms of $w_s$ and $\pi^a_s$  gives\[\max_{\pi\in\Delta^{|\S|}} r\cdot w[\pi] - \sum_{s\in\S} w[\pi]_s \sum_{as} \pi^a_s \log \pi^a_s = \max_{\pi^a_s}\left( r-h(\pi^a) \right) \cdot w[\pi]\] Writing out $w[\pi]$ explicitly in gives the regularized policy gradient
%which is exactly the objective function of the policy gradient  \eqref{eqn:lr_gdgd}.

%=======
\subsection{Equivalences}\hfill

{\bf Dual problem and policy gradient.} The dual problem \eqref{eqn:1r_d} is equivalent to the
policy gradient method. Let us parameterize $\mu^a_s = w_s \pi^a_s$ with $w_s =
\sum_{a\in\A}\mu^a_s$ and $\pi\in\Delta^{|\S|}$. Then the constraint in \eqref{eqn:1r_d} becomes
\[
1-\sum_{s\in\S} w_s = 0, \quad (I - (P^\pi)^\T)w = 0.
\]
which indicates that $w$ is the stationary distribution of $P^\pi$. After denoting this solution by
$w^\pi$ and plugging it into the objective function in \eqref{eqn:1r_d}, we transform the dual
problem to
\begin{equation}\label{eqn:1r_pg}
  \max_{\pi\in\Delta^{|\S|}} \left( r^\pi-h^\pi \right)^\T w^\pi.
\end{equation}
We can also view $r^\pi-h^\pi$ as a regularized reward by subtracting the entropy function $h^\pi$.
The average-reward under the new reward becomes
\[
\rho^\pi = \lim_{T\to\infty} \E\left[\frac1T\sum_{m=1}^T(r^{a_m}_{s_m} - h(\pi_{s_m}))\right].
\]
In this way, \eqref{eqn:1r_pg} is exactly the policy gradient method of this undiscounted
regularized MDP.

%% is
%% \begin{equation}\label{eqn:lr_gdgd}
%%   \trho = \max_{\pi\in\Delta^{|\S|}} \sum_{s\in\S}w^\pi_s \left(r^\pi_s - h^\pi_s\right)
%%   =\max_{\pi\in\Delta^{|\S|}} w^\pi  \cdot (r^\pi - h^\pi),
%% \end{equation}
%% which is equivalent to the dual problem \eqref{eqn:1r_d}.

{\bf Primal problem and Bellman equation.} The regularized average-reward Bellman equation for
$v,\rho$ is
\begin{equation}\label{eqn:1r_be}
  v = \max_{\pi\in \Delta^{|\S|}} r^\pi - \rho{\bf 1} +  P^\pi v - h^\pi.
\end{equation} 
In each component,
\[
v_s = \max_{\pi_s\in \Delta} \sum_{a\in\A}(r^a_s - \rho + \sum_{t\in\S}P^a_{st} v_t))\pi^a_s - h(\pi_s).
\]
By the Gibbs variational principle, the RHS is equal to $\log(\exp(\sum_{a\in\A}(r^a_s - \rho +
\sum_{t\in\S}P^a_{st} v_t)))$. Therefore, the regularized average-reward Bellman equation
can be written as the log-sum-exp form for $v,\rho$
\begin{equation}\label{eqn:lr_logsum}
  v_s = \log\left(\sum_{a\in\A} \exp\left(r^a_s - \rho +\sum_{t\in\S} P^a_{st}v_t\right)\right).
\end{equation}
Next, we show that the primal problem \eqref{eqn:1r_p0} is equivalent to solving \eqref{eqn:1r_be}.

{\em Bellman equation to primal problem.}  Let $\tv,\trho$ be the solution to the Bellman equation.
For $\tv,\trho$, there exists $\tpi$ s.t.,
\[
    \tv = r^{\tpi} - \trho{\bf 1} + P^{\tpi}\tv - h^{\tpi}.
\]
Then, for any $v$ satisfying the above constraints, the following inequality holds,
\[
    v \geq r^{\tpi} - \rho{\bf 1} + P^{\tpi}v - h^{\tpi}.
\]
Subtracting these two equations gives
\[
v-\tv \geq P^{\tpi}(v - \tv) - (\rho - \trho){\bf 1}
\]
Again let $w^*$ be the stationary distribution induced by the policy $\tpi$, then
$(w^*)^\T=(w^*)^\T P^{\tpi}$. Then multiplying $(w^*)^\T$ to the last equation yields
\[
(w^*)^\T(\rho - \trho){\bf 1} \geq 0
\]
Since we assume the MDP is unichain, the stationary distribution for any policy is strictly positive. This implies that
\[
\rho  \geq \trho
\]
holds for all $v$ satisfying the constraints in \eqref{eqn:1r_p0}. This proves that $\tv$ is the
minimizer of the primal problem \eqref{eqn:gr_p0}.

{\em Primal problem to Bellman equation.} Let $\tv, \trho = \argmin_{v,\rho} \rho$ be the minimizer
of the primal problem \eqref{eqn:1r_p0}. Besides, there exists a policy $\tpi$ s.t. $\trho =
(w^{\tpi})^\T(r^{\tpi} - h^{\tpi}).$ This can be seen from multiplying $(w^\pi)^\T$ to the
constraint $\tv \geq r^\pi - \trho {\bf 1} + P^\pi \tv -h^\pi $. Due to $(w^\pi)^\T{\bf 1} = 1,
(w^\pi)^\T P^\pi = (w^\pi)^\T$, one has
\[
\trho \geq (w^\pi)^\T(r^\pi - h^\pi)
\]
for all $\pi$. Let $\tpi = \argmin_\pi (w^\pi)^\T(r^\pi - h^\pi)$, then the minimizer $\trho =
(w^{\tpi})^\T(r^{\tpi} - h^{\tpi})$.

We now show that $(\tv,\trho)$ is also the solution to the Bellman equation \eqref{eqn:1r_be}
by contradiction. Assume $(\tv,\trho)$ does not satisfy \eqref{eqn:1r_be}, then there must exist $\bs$,
s.t., for $\forall \pi$
\begin{equation}\label{eqn: 1r_equiv3}
    \tv_{\bs} \geq (r^\pi - \trho +  P^\pi \tv - h^\pi)_{\bs} + \delta
\end{equation}
with some positive constant $\delta>0$. %Let $\tpi$ be the policy that has average-reward $\trho$, i.e., $\trho = (r^{\tpi})^\T w^{\tpi}$, where $w^{\tpi}$ is the stationary distribution satisfying $(w^{\tpi})^\T = (w^{\tpi})^\T P^{\tpi}$. 
The inequality \eqref{eqn: 1r_equiv3} also holds for $\tpi$, so one can write it in vector form,
\[
\tv \geq r^{\tpi} - \trho {\bf 1} +  P^{\tpi} \tv - h^{\tpi} + \tilde{\delta},
\]
where $\tilde{\delta}$ is a vector with $\delta$ on its $\bs$-th element and $0$ on all other elements. 
Then multiplying $(w^{\tpi})^\T$ to the above equation yields,
\[
 0 \geq  \delta w^{\tpi}_{\bs}.
\]
Note that the RHS is always $> 0$ because the stationary distribution $w^{\tpi}$ is strictly
positive, which leads to a contradiction. Hence we conclude that $\tv, \trho$ is also the solution
to the Bellman equation \eqref{eqn:1r_be}.

%The value function under the new reward becomes $v^\pi_s = \E[\sum_{m\geq 0} (r_{s_m} - h(\pi_{s_m}^a) - \rho^\pi) | s_0 = s]$, then the value function under the optimal policy satisfies the regularized average-reward optimal Bellman equaiton: \begin{equation}\label{eqn:lr_be} v_s^* = \max_{\pi\in \Delta^{|\S|}} r_s -\trho + \g \sum_{t\in\S} P^\pi_{st}\tv_t - h(\pi_s),\end{equation} with $\sum_{s\in\S} \tv_sw^*_s = 1$, which is equivalent to the primal problem \eqref{eqn:1r_p0}. The equivalence could be derived similarly as for the discounted regularized MDP.  The optimal value function could also be written as the log-sum-exp form, \begin{equation}\label{eqn:lr_logsum} v_s^*  = \log\left(\sum_{a\in\A} \exp(r_s - \trho +\g\sum_{t\in\S} P^a_{st}v_t^*)\right). \end{equation} which is equivalent to another primal problem \eqref{eqn:1r_p}.

%================================================
%\section{Conclusion}
%\LY{...}

%Some references:  \cite{abbasi2014linear,wang2020randomized}.

\bibliographystyle{abbrv}


\end{document}